\documentclass[b4,11pt,leqno]{article}

\usepackage{amsmath} 
\usepackage{amssymb} 

\usepackage{latexsym}

\usepackage[OT2,T1]{fontenc}

\numberwithin{equation}{section}

\numberwithin{equation}{section}  

\allowdisplaybreaks

\title{An overview of the history of  
\\ 
projective representations 
\\
(spin representations) of groups\footnote{This is a self English translation, with small corrections, of Appendix A of my book {\it Introduction to the theory of projective representations of groups}\, in Japanese \cite{[Hir2018]}.}
}

\author{Takeshi HIRAI}

\date*{}

\begin{document}

\pagenumbering{arabic}

\maketitle

\setcounter{section}{0}
\setcounter{page}{1}

{\small 

\tableofcontents

}

\vskip1.5em
To overview the history of spin representations (projective representations) of groups, here we choose the style to insert necessary explanations in a chronology in such a way that first we put age and main items in bold face, and cite main histotical papers, then, explain mainly its contents and important meaning. (More detailed overview of the history, containing the views from theoretical physics, was given in [HHoH].) 
Hoping to give the original flaver as far as possible, we quote original sentances of papers frequently even in German and in French. You can translate them by using website dictionaries if necessary.

\addcontentsline{toc}{section}{A.1\; \;Prehistory}

\setcounter{section}{1}
\setcounter{subsection}{1}
\vskip1.5em 
\noindent 
{\Large\bf A.1\; \;Prehistory}
\vskip1em 
\addcontentsline{toc}{subsection}{A.1.1\; \;Benjamin Olinde Rodrigues}

{\bf\large A.1.1\; \;Benjamin Olinde Rodrigues\footnote{Benjamin Olinde Rodrigues (6 October 1795\;--\;17 December 1851)}} 
\vskip.5em
{\bf A.1.1.1\; Rodrigues's paper}\; 
\vskip.2em
{\bf 1840}\; 
{\bf (Virtual discovery of quaternion, and Rodrigues expression of 
 spatial motions)} 

\vskip.2em
{\small 
[Rod] 
Olinde Rodrigues, {\it 
Des lois g\'eom\'etriques qui r\'egissent les d\'eplacements d'un syst\`eme solide dans l'espace, et la variation des coordonn\'ees provenant de ses d\'eplacements consid\'er\'es ind\'ependamment des causes qui peuvent les produire,} 
 Journal\;de Math\'ematiques\;Pures et Appliqu\'ees, {\bf 5}(1840), 380--440.
 }
 
\vskip.3em

[Translation of the title] Geometric rules which govern displacements of solid bodies in the space, and the changes of the coordinates coming from these displacements, considered independent of the cause of their generations. 
 \vskip1em

{\bf (1)}\; This long paper (61 pages) studies motions in the 3-dimensional Euclidean space $E^3$ systematically, starting from a system of axioms which imitates Euclid's Elements, so to say. Its highlight is to express the product of two spatial motions by calculation formulas obtained with the help of triangles on the sphere and expressed explicitly by means of trigonometric functions. This calculation formulas give substantially the rule of the product in the {\it quaternion algebra}.  
\vskip1em

{\bf (2)}\; Another important contribution is {\it Rodrigues expression}\, of a spatial rotation. This is understood, in recent years, more useful in practice than usual Euler angles. The reason why is 

$\bullet$ (apart from Euler angles) it has no gaps (jumps) of parameters and accordingly differential equations can be written with it easily, and 
 multiple rotations beyond $2\pi$ (360$^o$) can be continuously written without problems, and 

$\bullet$  algorithms with less computational complexity can be given, etc.

\vskip1em

{\bf (3)}\; Aside from these things, there appear here substantially the double covering of the rotation group $SO(3)$ and its 2-dimensional spin representation (doubly-valued representation) as is explained below.

\vskip1.2em

{\bf A.1.1.2\; More detailed explanation.}\; 

The original description of Rodrigues himself is mainly geometrical by using triangles on a sphere. But, here we explain the contents of a paper written about 180 years ago, in the world of modern mathematics, so please permit us to use, for convenience, the quaternion itself and modern mathematical terminology. 

Let the standard units of quaternion be $1, {\boldsymbol i},{\boldsymbol j},{\boldsymbol k}$ (see (\ref{2013-01-14-1}) below). The quatenion is given as ${\boldsymbol H}={\boldsymbol R} 1+{\boldsymbol R}\,{\boldsymbol i}+{\boldsymbol R} {\boldsymbol j}+{\boldsymbol R} {\boldsymbol k}$，and the total set of pure quaternion numbers as ${\boldsymbol H}_-={\boldsymbol R} {\boldsymbol i}+{\boldsymbol R} {\boldsymbol j}+{\boldsymbol R} {\boldsymbol k}$．
The length of ${\boldsymbol x} =x_0+x_1{\boldsymbol i}+x_2{\boldsymbol j}+x_3{\boldsymbol k}\;(x_j\in{\boldsymbol R})$ is defined as $\|{\boldsymbol x}\|:=\sqrt{x_0^{\,2}+x_1^{\,2}+x_2^{\,2}+x_3^{\,2}}$, and the set of all quaternion numbers with length 1 is denoted by ${\boldsymbol B}$, and we put ${\boldsymbol B}_-:={\boldsymbol H}_-\cap {\boldsymbol B}$. Any element of ${\boldsymbol H}_-$ is expressed as $\phi{\boldsymbol w}\;(\phi\in{\boldsymbol R},\;{\boldsymbol w}\in{\boldsymbol B}_-)$. 

The following lemma gives elementary properties necessary in the following. 

\vskip1.2em

{\bf Lemma A.1.1.}  {\it Put\; 
$\overline{{\boldsymbol x}}:=x_0-x_1{\boldsymbol i}-x_2{\boldsymbol j}-x_3{\boldsymbol k}$. 

{\rm (i)}\;\;\; 
${\boldsymbol x}\,\overline{{\boldsymbol x}}=\|{\boldsymbol x}\|^2$\; and, for ${\boldsymbol x}\ne 0$, we have ${\boldsymbol x}^{-1}=\overline{{\boldsymbol x}}/\|{\boldsymbol x}\|^2$. 

{\rm (ii)}\;\; $\overline{{\boldsymbol x}{\boldsymbol y}}=\overline{{\boldsymbol y}}\cdot\overline{{\boldsymbol x}}.$ 

{\rm (iii)}\; $\|{\boldsymbol x}{\boldsymbol y}\|=\|{\boldsymbol x}\|\cdot\|{\boldsymbol y}\|\;\;({\boldsymbol x},{\boldsymbol y}\in{\boldsymbol H}).$
}
\vskip1em

{\it Proof.}\; 
By using (i) and (ii), (iii) is given as 
\vskip.5em 
\qquad
 $\|{\boldsymbol x}{\boldsymbol y}\|^2=({\boldsymbol x}{\boldsymbol y})\,\overline{{\boldsymbol x}{\boldsymbol y}}= 
({\boldsymbol x}{\boldsymbol y})(\overline{{\boldsymbol y}}\cdot\overline{{\boldsymbol x}})= 
 \|{\boldsymbol x}\|^2\cdot\|{\boldsymbol y}\|^2.$ 
 \hfill 
 $\Box$\;
 \vskip1.2em
 
This lemma shows that ${\boldsymbol H}$ is a skew field, and the map \,${\boldsymbol x}\mapsto \overline{{\boldsymbol x}}$\, reverses the order of the product.  By the assertion (iii), ${\boldsymbol B}$ is closed with respect to the product, and since $\|u\|=1$ gives $\|u^{-1}\|=1$, we know that ${\boldsymbol B}$ is a group.

 \vskip1.2em

{\bf Lemma A.1.2.\footnote{This is a corrected version of Lemma A.1.2 in \cite{[Hir2018]}.}}\; 
{\it If we consider  
${\boldsymbol H}$ as an algebra over ${\boldsymbol R}$, then there exists an isomorphism from its complexification ${\boldsymbol H}_{\boldsymbol C}=H\otimes_{\boldsymbol R} {\boldsymbol C}$ onto the full matrix algebra $M(2,{\boldsymbol C})$ of complex matrices of order 2. An isomorphism $\Phi$ is given by the correspondence: 
\begin{eqnarray*}
{\boldsymbol i}\to I= 
{\small 
\begin{pmatrix}
0 & -1 \\ 1 &0 
\end{pmatrix}
},
\quad 
{\boldsymbol j}\to J=
{\small 
\begin{pmatrix}
0&i \\ i&0
\end{pmatrix}
},
\quad 
 {\boldsymbol k}\to K=
{\small 
\begin{pmatrix}
-i&0 \\ 0&i
\end{pmatrix}}
.
\end{eqnarray*}
}
\vskip.2em

{\it Proof.}\; 
By simple calculations, it is shown that the matrices $I,J,K$ satisfy the equations corresponding to the so-called {\it fundamental formula}\; ${\boldsymbol i}^2={\boldsymbol j}^2={\boldsymbol k}^2=-1,\;{\boldsymbol i}{\boldsymbol j}{\boldsymbol k}=-1$. Together with the unit matrix $E_2$, they form a basis over ${\boldsymbol C}$ of $M(2,{\boldsymbol C})$. 
\hfill 
$\Box$\;

 \vskip1em

{\bf A.1.1.3\; Rodrigues expression of rotations.} 

We identify the 3-dimensional Euclidean sapce $E^3$ with ${\boldsymbol H}_-$ under the correspondence ${\boldsymbol H}_-\ni{\boldsymbol x} =x_1{\boldsymbol i}+x_2{\boldsymbol j}+x_3{\boldsymbol k}\longleftrightarrow x={}^t(x_1,x_2,x_3)\in E^3\;(x_j\in{\boldsymbol R})$, where a point in $E^3$ is denoted with a column vector $x$. 
Since it can be proved that a rotation $R$ arround the origin ${\bf 0}\in {\boldsymbol H}_-$ has necessarily an axis of rotation,\footnote{Any rotation in odd-dimensional space necessarily has an axix of rotation.} it can be expressed by means of the axis of rotation ${\boldsymbol w}\in{\boldsymbol B}_-$ and the angle $\phi\in{\boldsymbol R}$ measured as a right-handed screw rotation. Here we accept any angle, not restricted as $|\phi|\le 2\pi$. Take the vector $\phi{\boldsymbol w}\in{\boldsymbol H}_-$ as a parameter of $R$, and denote $R$ as  
\begin{eqnarray} 
\label{2014-09-07-1}
R=R(\phi{\boldsymbol w})\qquad(\phi{\boldsymbol w}\in{\boldsymbol H}_-). 
\end{eqnarray} 
This is called as {\it Rodrigues expression}\, of rotation. 

Now let us realize the rotation $R(\phi{\boldsymbol w})$ explicitly by means of $\phi$ and ${\boldsymbol w}$.

 \vskip1.2em

{\bf Lemma A.1.3.} {\it 
The unit ball ${\boldsymbol B}$ in 
${\boldsymbol H}$ is a group with respect to the product, and it gives a universal covering group of the rotation group $SO(3)$. 
}
 \vskip.8em
 
{\it Proof.}\; For a ${\boldsymbol u}\in{\boldsymbol B}$, define a map
\begin{eqnarray*}
\label{2013-01-15-1}
\nonumber
\Psi({\boldsymbol u}):\, {\boldsymbol H}_-\ni{\boldsymbol x}\;\mapsto\; {\boldsymbol x}'={\boldsymbol u}{\boldsymbol x}{\boldsymbol u}^{-1}\in{\boldsymbol H}_-.
\end{eqnarray*}
Then, by Lemma A.1.1, $\|{\boldsymbol x}'\|=\|{\boldsymbol u}\|\cdot\|{\boldsymbol x}\|\cdot\|{\boldsymbol u}^{-1}\|=\|{\boldsymbol x}\|$. Hence the corresponding map  
 $$
 \Psi'({\boldsymbol u}):\,E^3\ni x\longmapsto x'=gx\in E^3
 $$ 
conserves the length, that is, $g=g({\boldsymbol u})$ belongs to the length preserving group $O(3)$ of $E^3\cong {\boldsymbol H}_-$. Moreover, $g$ can be connected to the unit matrix $E_3$ continuously, and so it belongs to the connected component $SO(3)$ containing $E_3$, or it gives a rotation in $E^3$ arround the origin. 
The kernel of the homomorphism $\Psi': {\boldsymbol B}\ni{\boldsymbol u}\mapsto g({\boldsymbol u})\in SO(3)$ is ${\Psi'}^{\,-1}(\{E_3\})=\{\pm 1\}$, and so ${\boldsymbol B}/\{\pm 1\}\cong SO(3)$. Moreover, the 3-dimensional ball is simply connected, and so ${\boldsymbol B}$ is simply connected. Accordingly ${\boldsymbol B}$ is a universal covering group of $SO(3)$. 
 \hfill
 $\Box$\;\; 

  \vskip1.2em

{\bf Lemma A.1.4.}   
{\rm (i)}\; 
{\it For any ${\boldsymbol x}\in{\boldsymbol H}_-$, there exists a ${\boldsymbol u}\in {\boldsymbol B}$ such that ${\boldsymbol x}={\boldsymbol u}\cdot \theta{\boldsymbol i}\cdot {\boldsymbol u}^{-1}$ with $\theta=\|{\boldsymbol x}\|$. 

{\rm (ii)} The subset ${\boldsymbol B}_-={\boldsymbol B}\cap{\boldsymbol H}_-$ of ${\boldsymbol B}$ is a conjugacy class of ${\boldsymbol B}$ and 
${\boldsymbol B}_-=\{{\boldsymbol x}\in {\boldsymbol B}\,;\,{\boldsymbol x}^2=-1\}$. 
}
 \vskip1em
 
{\it Proof.}\;  
(i)\; As is well known, any element $x\in E^3$ is sent to ${}^t(\theta,0,0)$ with $\theta=\|x\|$ under the action of $SO(3)$. By the correspondence 
 ${\boldsymbol H}_-\cong E^3$ and $\Psi':\,{\boldsymbol B}\ni {\boldsymbol u} \mapsto g({\boldsymbol u})\in SO(3),$ in the proof of the previous lemma, pull back this to the world of ${\boldsymbol H}_-$ and ${\boldsymbol B}$. 

(ii) Apply (i) to an ${\boldsymbol x}\in{\boldsymbol B}_-$. Then $\theta=\|{\boldsymbol x}\|=1$, and the element in ${\boldsymbol H}_-$ corresponding to ${}^t(1,0,0)\in E^3$ is ${\boldsymbol i}$. Hence ${\boldsymbol x}={\boldsymbol u}{\boldsymbol i}{\boldsymbol u}^{-1}\;(\exists {\boldsymbol u}\in B)$. 
Also express ${\boldsymbol x}\in {\boldsymbol B}$ as ${\boldsymbol x}=\alpha+\beta{\boldsymbol i}+\gamma{\boldsymbol j}+\delta{\boldsymbol k}\;(\alpha,\beta,\gamma,\delta\in{\boldsymbol R})$, then $\alpha^2+\beta^2+\gamma^2+\delta^2=1$. Apply here the condition ${\boldsymbol x}^2=-1$, then we see that it is equivalent to the condition $\alpha=0$. 
 \hfill
 $\Box$\;\;

\vskip1.2em

{\bf Lemma A.1.5.}  
{\it For\; 
${\boldsymbol x}\in{\boldsymbol H}$, the infinite series 
$$ 
\exp({\boldsymbol x}):=
\sum_{0\leqslant k< \infty}\frac{1}{k!}\,{\boldsymbol x}^k
$$  
is absolutely convergent, and for 
${\boldsymbol v}\in{\boldsymbol H}_-$, we have $\exp({\boldsymbol v})\in{\boldsymbol B}$\, or $\|\exp({\boldsymbol v})\|=1$. Moreover expresse ${\boldsymbol v}$ as ${\boldsymbol v}=\phi{\boldsymbol w}\;({\boldsymbol w}\in {\boldsymbol B}_-,\,\phi\in {\boldsymbol R})$, then \,
\vskip.5em 
\hspace*{4ex} $\exp(\phi{\boldsymbol w})= \cos\phi+\sin\phi\cdot{\boldsymbol w},\quad\exp\big((\phi+2k\pi) {\boldsymbol w})=\exp(\phi {\boldsymbol w})\;\;(k\in{\boldsymbol Z})$. 

}

\vskip1em

{\it Proof.}\;  From\; 
$\sum_{0\leqslant k<\infty}\|\frac{1}{k!}{\boldsymbol x}^k\|=\sum_{0\leqslant k< \infty}\frac{1}{k!}\|{\boldsymbol x}\|^k=\exp\big(\|{\boldsymbol x}\|)<\infty$,\, the absolute convergence follows.  
 Also, from the expression\; 
${\boldsymbol v}=\phi{\boldsymbol w}$, by using\;  ${\boldsymbol w}^2=-1$, we obtain  
\begin{eqnarray*}
\hspace*{8ex}
\exp(\phi{\boldsymbol w})=\sum_{p\geqslant 0}\frac{(-1)^p}{(2p)!}\phi^{2p} + \sum_{p\geqslant 0}\frac{(-1)^p}{(2p+1)!}\phi^{2p+1}{\boldsymbol w} 
=\cos\phi+\sin\phi\cdot{\boldsymbol w}.
\hspace{7ex}\Box
\end{eqnarray*}
\vskip.2em 

{\bf Theorem A.1.6.} (Rodrigues expression of a rotation)\; 
{\it 
Make the angle of rotation to a half as\;
 $\phi\to\tfrac{1}{2}\phi$\; and then take a composition map  \;$\Psi'\circ \exp:\, {\boldsymbol H}_-\to SO(3)$. This gives the rotation $R(\phi{\boldsymbol w})\in SO(3)$ as    
\begin{eqnarray*}
\label{2013-01-14-2} 
&&
{\boldsymbol H}_-\ni\phi{\boldsymbol w} \;\mapsto\; \Psi'\big(\exp(\tfrac{1}{2}\phi{\boldsymbol w})\big)=\Psi'\big(\cos(\tfrac{1}{2}\phi)+\sin(\tfrac{1}{2}\phi){\boldsymbol w}\big)\in SO(3).
\end{eqnarray*}
Or the right hand side expresses the element in the rotation group  
 $SO(3)$ with rotation axis ${\boldsymbol w}$ and rotation angle $\phi$.
} 

\vskip.8em

{\it Proof.}\; It is sufficient to prove the assertion in the case of  
${\boldsymbol w}={\boldsymbol i}$ (Why ?). Denote by $\rho$ the rotation $\Psi'\big(\exp(\tfrac{1}{2}\phi\,{\boldsymbol i})\big)$. Then, $\rho({\boldsymbol i})={\boldsymbol i}$. Thus the rotation axis of $\rho$ is ${\boldsymbol i}$, and so $\rho$ is expressed as a rotation in the  $({\boldsymbol j},{\boldsymbol k})$-plane perpendicular to ${\boldsymbol i}$. By using the double angle formula for trigonometric function, we obtain 
\begin{eqnarray}
\label{2013-01-14-3}
\nonumber
\big(\rho({\boldsymbol j}),\rho({\boldsymbol k})\big)= ({\boldsymbol j},{\boldsymbol k})
{\small 
\begin{pmatrix}
\cos \phi &-\sin\phi 
 \\ 
 \sin\phi & \;\;\cos\phi
\end{pmatrix}
}.
\end{eqnarray}
Hence $\rho$ is a rotation in the $({\boldsymbol j},{\boldsymbol k})$-plane with angle $\phi$.
 \hfill
 $\Box$\;\; 
\vskip1.2em

{\bf Important remark.} 

Under the map $\Psi'$, the angle  $\theta=\tfrac{1}{2}\phi$ in the parameter space ${\boldsymbol H}_-$ appears twice as  $2\theta=\phi$ for the rotation angle in $E^3$. This means that $\Psi':\,{\boldsymbol B}\ni {\boldsymbol u} \mapsto g=\Psi'({\boldsymbol u})\in SO(3)$ is a map {\bf 2\,:\,1} and the group ${\boldsymbol B}$ is a double covering group of $SO(3)$.
Also, reverse the perspective and consider the map (locally univalent but globally doubly-valued) as 
\begin{eqnarray}
\label{2019-08-01-1}
\pi:\;SO(3)\ni g=g({\boldsymbol u}) \,\longmapsto\, {\boldsymbol u} \in{\boldsymbol B}\subset {\boldsymbol H}.
\end{eqnarray} 
Connect with this the isomorphism $\Phi:{\boldsymbol H}\otimes_{\boldsymbol R}{\boldsymbol C}\to M(2,{\boldsymbol C})$, then there appears the 2-dimensional projective representation (the so-called {\it spin representation})\; 
$\Phi\circ \pi$\, of the rotation group $SO(3)$.

\newpage

{\bf A.1.1.4. Product of two rotations and calculation rule in quaternion.} 

It is very natural from the above dicussions to imagine that, when the product 
$$
R(\phi{\boldsymbol w})R(\phi'{\boldsymbol w}')\qquad({\boldsymbol w},{\boldsymbol w}'\in{\boldsymbol B}_-,\;\phi,\phi'\in{\boldsymbol R})
$$ 
of two rotations is expressed as $R(\phi''{\boldsymbol w}'')$, the rule of calculating $(\phi'',{\boldsymbol w}'')$ from $(\phi,{\boldsymbol w})$ and $(\phi',{\boldsymbol w}')$ comes out from the product rule in quaternion. (Note that, to modernize terminologies and summarize the total story, we used quaternion from the beginning of the above discussons. We sincerely apporogise that the historical time goes back and forth as this.) 

In the above paper of Rodrigues, the word \lq\lq quaternion'' never appeared, but only geometric calculations concerning triangles on a sphere. 
However the calculation formula for the composition of two rotations are explicitly written out. This calculation rule is equivalent to the following formula in the quaterinion ${\boldsymbol H}$\,: 
\begin{equation}
\label{2008-02-20-33}
{\textstyle 
\big(\cos({\frac{1}{2}}\phi)+\sin(\frac{1}{2}\phi){\boldsymbol w}\big)
\big(\cos(\frac{1}{2}\phi')+\sin(\frac{1}{2}\phi'){\boldsymbol w}'\big)
=
\cos(\frac{1}{2}\phi'')+\sin(\frac{1}{2}\phi''){\boldsymbol w}''.
}
\end{equation}

\vskip.2em

{\bf Note.} Hoping to get a formula which is completely univalent, if we put some restriction to the range of the variable $\phi$ in $\tfrac{1}{2}\phi {\boldsymbol w}$, then some breaks of parameter $\phi$ will appear. Accordingly it is better and natural to endure multivalency and take ${\boldsymbol H}_{\!-}$ as the parameter space, and admit multiple rotations arround a rotation axis ${\boldsymbol w}$ (Cf. [H5]).

\vskip2em

\addcontentsline{toc}{subsection}{A.1.2\; \;William Rowan Hamilton} 

{\bf\large A.1.2.\; \;William Rowan Hamilton\footnote{William Rowan Hamilton (4 August 1805\;--\;2 September 1865)}} 
\vskip.5em 

{\bf 1843} {\bf （Discovery of quaternion）}

{\small
[Ham]  
W.R.\;Hamilton, {\it  
On a new species of imaginary quantities connected with the theory of quaternions,} Proc. Royal Irish Acad., {\bf 2}(1843), 424--434.
} 
\vskip.8em

Starting from a certain point of his research life, 
Hamilton continued long time, as his research objective, to extend complex field  \,${\boldsymbol C}={\boldsymbol R} 1+{\boldsymbol R} i,\;i=\sqrt{-1}$, to a vector space and make it possible to divide a vector by another vector (he called this as a division by a vector). As the first step, he tried the case of adding two imaginary units to the real field ${\boldsymbol R}$, and after a considerable effort, he could succeeded finally, on the 16th of October 1843, to discover that, by adding three imaginary units, ${\boldsymbol i}, {\boldsymbol j}, {\boldsymbol k}$, and introducing the following calculation rule, which Hamilton called {\it fundamental formula}, it is possible\,:     
\begin{eqnarray}
\label{2013-01-14-1}
\qquad
{\boldsymbol i}^2={\boldsymbol j}^2={\boldsymbol k}^2=-1,\quad {\boldsymbol i}{\boldsymbol j}{\boldsymbol k}=-1\qquad\mbox{\rm ({\rm fundamental formula}).}
\end{eqnarray}

\vskip.3em

After such a big discovery, firmly believing that \lq\lq {\it the quaternion should be a great discovery of the century}\,'', Hamilton was making efforts to show usefulness of the quaternion and to expand the theory to the scientific world. However, he could not succeeded to express a 3-dimensional rotation completely well 
by means of the quaternion. The main reason of this failure is that, with the pride and the stubbornness specific to the great dicoverer himself, he tried to find out a  {\bf 1\,:\,1} correspondence between the unit ball ${\boldsymbol B}\cong S^3$ of quaternion numbers and the total set $SO(3)$ of 3-dimensional rotations. 

Actually, in the case of complex number field and 2-dimensional rotations, under the natural correspondence ${\boldsymbol C}={\boldsymbol R} +{\boldsymbol R} i\ni x+yi\longleftrightarrow {}^t(x,y)\in E^2$, identify ${\boldsymbol C}$ with the 2-dimensional Eucledean space $E^2$, then an element $w=e^{i\theta},\,0\le \theta<2\pi,$ of ${\boldsymbol B}_2:=\{w\in{\boldsymbol C}\,;\,|w|=1\}$ naturally corresponds 1\,:\,1 way to a 2-dimensional rotation of angle $\theta$ under the simple multiplication as\;  
$z=x+yi\mapsto z'=e^{i\theta}z=x'+y'i$, where  
$$
\begin{pmatrix}
x' \\ y'
\end{pmatrix}
= 
u_2(\theta)
\begin{pmatrix}
x \\ y
\end{pmatrix},\quad u_2(\theta)= 
\begin{pmatrix}
\cos \theta& -\sin\theta 
 \\ 
 \sin\theta &\;\;\cos\theta
\end{pmatrix}.
$$
Or, in the case of 2-dimensional rotations the natural isomorphism ${\boldsymbol B}_2 \ni e^{i\theta}\to u_2(\theta)\in SO(2)$ is realized. 

\vskip1.2em

{\bf Remark A.1.2.1. Comentary from a point of view of the history of science.}

We quote here a paper explaning a rather complicated situation from a point of view of the history of science arround Rodrigues, Hamilton and the quaternion. From this paper, I myself got a lot of new important knowledge which I have never seen. 
\vskip.7em 

[Alt]\, 
S.L.\;Altmann, 
{\it Hamilton, Rodrigues, and the quaternion scandal, 
\, 
{\small 
What went wrong with one of the major mathematical discoveries of the nineteenth century},
} Mathematics Magazine, {\bf 62}(1989), 291--308.

\vskip.7em

There is a reality, which cannot be put aside under the simple word \lq{\it Irony of the History}, that the virtual discovery of quaternion of Rodrigues and his other important mathematical works were not taken into account and ignored totally. 

Historically, as a Jewish French, he suffered by the anti-Semitism of Cathoric Charch Power after Restoration of Imperial Rule something like \lq{\it Exclude Jews from public institutions}\,', and so he could not find any teaching position and was obliged to be a banker as his family business. But at his age over 40, he could publish the above excellent, long mathematical paper [Rod], with author's name Olinde Rodrigues, in Journal de Math\'ematiques Pures et Appliqu\'ees.  

\'Elie Cartan cited this paper in his book [Car2] on the theory of spinor ({\it spineur}\, in French), on the subject of 2-dimensional spin representation of $SO(3)$ or of $\mathfrak{s}\mathfrak{o}(3)$, but he misunderstood that the paper is written by two coauthors called Olinde and Rodrigues respectively. In reality, the first name of Rodrigues is  
Benjamin and the middle name is Olinde, which is added in his boyhood by his father after the official order of Cathoric Church Power that \lq{\it Everyone should have middle name}', but this middle name {\it Olinde}\, does not follow the cathoric style and is very unusual (one of my French friends pointed out it to me). It may be a reason why \'{E}. Cartan mistook as two coauthors. Rodrigues himself used his middle name Olinde as the first name in his mathmatical papers. 

Before and after this paper, between 1838 and 1843, Rodrigues published 7 papers altogether in the same Journal. They contain short papers but every paper is good one with sharp cuts, I feel. He pulished them in between 43 and 48 ages, with the long silence over 20 years after his university life.

Why\,? I don't know, but these science-historical facts have been totally ignored until recently. However the mathematician Benjamin Olinde Rodrigues (1795--1851) now becomes known and his honor is recovered (Cf. [AlOr]). Myself I studied him in many ways and made a report on B.\,O.\,Rodrigues in the 22nd symposium on the history of mathematics, 2011, held at Tsuda University (Cf. [H5]).

\vskip2.5em

{\Large\bf A.2.\; \;History}

\addcontentsline{toc}{section}{A.2\; \;History}

\setcounter{section}{2}
\setcounter{subsection}{0}
\vskip.5em 

\addcontentsline{toc}{subsection}{A.2.1\; \;Creation of theory of representations of groups}

\vskip.8em 

{\large\bf A.2.1\; \;Creation of theory of representations of groups}

\vskip.5em
{\bf 1896\; (It begins with the theory of characters and group determinants)}

\vskip.4em

{\small 
[Fro1]\;    
F.\;Frobenius, {\it 
\"Uber Gruppencharaktere,} 
Sitzungsberichte der K\"oniglich Preu\ss ischen Akademie der Wissenschaften zu Berlin, {\bf 1896}, pp.985--1021.

\vskip.3ex

[Fro2]\;  
---, {\it 
\"Uber die Primfactoren der Gruppendeterminante,} 
ibid., {\bf 1896}, pp.1343--1382.

\vskip.3ex

[Fro3]\; 
---, {\it 
\"Uber die Darstellung der endlichen Gruppen durch lineare Substitutionen,} 
ibid., {\bf 1897}, pp.944--1015.

} 

\vskip.5em

Dedekind was studying characters of non-abelian finite groups $G$, starting from 1880, and put some questions to Frobenius.\footnote{Ferdinand Georg Frobenius (26 October 1849\;--\;3 August 1917)}  From this occasin Frobenius has started his study on   the theory of characters of finite non-abelian groups and so on. The first results were published in papers [Fro1,\,1896], [Fro2,\,1896]. These are, above all, the beginning of the systematic studies of the theory of linear representations (of finite groups). Just after, in the next paper [Fro3,\,1897], the study of linear representations started and the relations to group characters and group determinants were clarified.  
\vskip.8em

Frobenius defined in [Fro1] {\it characters}\, of non-abelian group by a system of equations given by means of orders of several kinds of conjugacy classes, and gave a complete set of its solutions. He studied in [Fro2] {\it group determinant}\, in a purely algebraic way. Next year, in [Fro3], he began to study actually {\it linear representations}\, of (finite) groups, and proved that the character defined in [Fro1] is equal to the trace character of irreducible linear repreentations, and the group determinant in [Fro2] comes out from the regular representation on $\ell^2(G)$ and rules over, from an algebraic point of view, its decomposition into irreducible representations. 

\vskip.8em

Afterwords, in 12 years until 1906, he wrote about 14 papers on the theory of characters and linear representations of groups (2 papers on 1906 are coauthored with I.\,Schur) which gave main principal results in the theory of group representations. This is about 12 years starting from 46 years of age or so (he was nominated as professor of University of Berlin in 1893 and became a menber of Prussian Royal Academy of Science). I read these papers in detail and reported on them in the symposium on the history of mathematcs at Tsuda University four times (Cf. [H1]). 

The structure of Frobenius' theory is very strict and algebraic, but very difficult to read. After him, Burnside\footnote{William Burnside (2 July 1852\;--\;21 August 1927)} gave in [Bur1,\,1898]\,--\,[Bur2,\,1898] another way of approach, and Schur\footnote{Issai Schur (10 January 1875\;--\;10 January 1941)}, a student of Frobenius, reconstruct the theory of linear representations and characters of a finite group in a modern style [Sch6,\,1905].   

\vskip1.5em

{\bf 1897}\; {\bf (Molien's work, overlapped partially with the work of Frobenius)}

\vskip.4em
{\small 
[Mol2]  
T.\;Molien, {\it 
Eine Bemerkung zur Theorie der homogenen Substitutionensgruppe,}
 Sitzungsberichte der Dorpater Naturforscher-Gesellschaft, {\bf 11}(1897), 259--274. 

}

\vskip.4em
Theodor Molien\footnote{Theodor Molien (10 September 1861\;--\;25 December 1941)} (Russian name, Fedor Eduardovich Molin) obtained his doctor's degree at Dorpat Imperial University, a frontier of German academic study of that time. The work of his thesis was published as [Mol1, 1892] and the main contribution is a proof of a part of the following Wedderburn's Theorem in the case of the coefficient field  $K={\boldsymbol C}$\,:
\begin{quotation}
{\it 
A simple algebra over a commutative field $K$ is isomorphic to the full matrix algebra $M(n,D)$ of a skew field $D$ over $K$. 
}
\end{quotation}

Later he got a teaching position in that university, he wrote a paper [Mol2,\,1897] in 
Bulletin of Dorpat Natural Science Researcher's Association, which contains an important result proved by using the above result in [Mol1], independent of Frobenius. (Actually Frobenius utilized the result of himself on associative algebra obtined erlier.) 
In \S 4 of [Fro3,\,1897], Frobenius noticed this fact refering papers [Mol1] and [Mol2]. 

A paper [Mol3, 1898] was published in Bulletin of Prussian Academy in Berlin under the introduction by Frobenius. It seems that, afterwards Frobenius looked for some teaching position for this young talented reseacher, but it did not succeed. Later, when Dorpat was included in Soviet Union, Molien chose to remain there and was sent to Siberia and there he wrote many textbooks in mathematics (Cf. a historical study [H4]).

\newpage 

\addcontentsline{toc}{subsection}{A.2.2\; \;Beginning of projective (or spin) representations of groups}

{\Large\bf \;A.2.2\; \;Beginning of projective (or spin) representations of groups}

\vskip.8em 

It was 1904, only 8 years later from the creation of the theory of linear repersentations of groups in 1896, that the theory of {\it projective representations}\, or {\it spin representations}\, of groups begun.  

\vskip.6em 

{\bf 1904} {\bf (Projective representations [representations by linear fractional transformations] or spin representations of groups)}
\vskip.4em

{\small 
[Sch1]\; J.\;Schur,\; {\it 
\"Uber die Darstellung der endlichen Gruppen durch 
  gebrochene lineare Substitutionen,}   J. f\"ur die reine und angewante 
  Mathematik, {\bf 127}(1904), 20--50.
 
[Title] On representations of finite groups through linear fractional transformations.
}

 \vskip.4em

After his teacher Frobenius opened wide the gate to the thoery of {\it linear representations}\, of groups, Schur opend the gate to the theory of {\it projective representations}\, or {\it spin representatrions}.  The former is a homomorphism into $GL(n,K),\,K={\boldsymbol C}$ or ${\boldsymbol R},$ and the latter is a homomorphism into $PGL(n,K)=GL(n,K)/K^\times$. But Schur himself noticed at the top of this paper as 
 
\begin{quotation}
Das Problem der Bestimmung aller endlichen Gruppen lineare Substitutionen bei gegebener Variabelnzahl $n\;(n>1)$ geh\"ort 
zu den schwierigsten Problemen der Algebra und hat bis jetzt nur f\"ur die bin\"aren und tern\"aren Substitutionsgruppen seine vollst\"andige L\"osung gefunden. F\"ur den allgemeinen Fall ist nur bekannt, da\ss\;die Anzahl der in Betracht kommenden Typen von Gruppen eine endliche ist; dagegen fehlt noch jede \"Ubersicht \"uber die charakteristischen Eigenschaften dieser Gruppen.

Die Umkehrung dieses Problems bildet in einem gewissen Sinne die Aufgabe: alle Gruppen von h\"ochstens $h$ ganzen oder gebrochenen linearen Substitutionen zu finden, die einer gegebenen endlichen Gruppe $\mathfrak{H}$ der Ordnung $h$ ein- order mehrstufig isomorph sind, oder auch, wie man sagt, alle Darstellungen der Gruppe $\mathfrak{H}$ durch lineare Substitutionen zu bestimmen. 
\end{quotation}

(Translation)\, 
The problem of determining all the finite groups of linear transformations of given degree $n\;(n>1)$ belongs to the most difficult problems in algebra. And a complete answer is obtained only for $n=2, 3$. In general cases, we know only that the number of types of such finite groups is finite, and don't have any prospect about characterizing properties of such groups.

If we think about it, this problem suggests the following research subject: 
\vskip.6em
{\it 
Find all groups consisting of at most $h$ number of linear transformations or 
of projective transformations, which is isomorphic or homomorphic 
 to a given finite group $\mathfrak{H}$ of order $h$. Or, in other words, determine all the linear or projective representations of\, $\mathfrak{H}$. 
}

\vskip1em

Therefore Schur himself has an idea to {\it determine all the finite groups contained in $GL(n,{\boldsymbol C})$}. So the way of stating the results in this paper is that he put his weight to know the image group $\pi(G)\subset GL(n,{\boldsymbol C})$ of representation $\pi$ of $G$, which is more stronger than the motivation to represent $G$ or, so to say, the representation $\pi$ itself as a map. Accordingly any symbols denoting representation itself (like $\pi$) do not used. 

\vskip1em

{\bf 1.\;Projective representation.}\; 
Here we give definition of projective representations, essentially same as that of Schur. 
A {\it projective representation}\, of a group $G$ is a map $\pi$ which gives, for $g\in G$, a linear transformation $\pi(g)$ over ${\boldsymbol C}$ in such a way that 
\begin{eqnarray}
\label{2013-01-18-1}
\pi(g)\pi(h)=r_{g,h}\,\pi(gh)\quad(g,h\in G,\; r_{g,h}\in{\boldsymbol C}^\times),
\end{eqnarray}
and the function $r_{g,h}$ on $G\times G$ is called the {\it factor set}\, of $\pi$.

On the other hand, a function $r$ on $G\times G$ with values in ${\boldsymbol C}^\times:=\{z\in{\boldsymbol C};\,z\ne 0\}$ satisfying 
\begin{eqnarray}
\label{2013-01-19-1}
r_{k,gh}\,r_{g,h}=r_{k,g}\,r_{kg,h}\quad(k,g,h\in G), 
\end{eqnarray}
is called as a {\it ${\boldsymbol C}^\times$-valued 2-cocycle on}\, $G$. 
The product of 2-cocycles is also a 2-cocycle, and the quotient of the set of all 2-cocycles modulo the equivalence relation $r_{g,h}\approx r'_{g,h}:=r_{g,h}\cdot (\lambda_g\lambda_h/\lambda_{gh}),$ with $\lambda_g\in{\boldsymbol C}^\times\;(g\in G)$, is denoted by $H^2(G,{\boldsymbol C}^\times)$ and is called {\it Schur multiplier}\, of $G$. 

\vskip1em

In the paper [Sch1], the fundmental theory of projective representations of finite groups are discussed. 
Nowadays, {\it a central extention}\, $G'$ of $G$ by a commutative group $Z$ is defined by an exact sequence  
\begin{eqnarray}
\label{2013-01-19-2}
1 \to Z\to G'\stackrel{\Phi\,}{\,\to\,}G\to 1\quad{\rm (exact)},
\end{eqnarray}
with a homomorphism $\Phi:\,G'\to G$, and $Z\hookrightarrow Z(G')$, the center of $G'$. 
Take in $G'$ a section ${\cal S}$ of $G$ as $G\ni g\to s(g)\in {\cal S}\subset G'$, then  
\begin{eqnarray}
\label{2013-01-18-2}
s(g)s(h)=z_{g,h}\,s(gh)\quad(g,h\in G,\;\exists z_{g,h}\in Z). 
\end{eqnarray}

Now, for a linear representation $\Pi$ of $G'$, put $\pi(g):=\Pi\big(s(g)\big)\;(g\in G)$. Then, for $g,h\in G$, 
\begin{eqnarray*}
\label{2013-01-18-3}
\pi(g)\pi(h)
\!\!&=&\!\!
\Pi\big(s(g)\big)\Pi\big(s(h)\big)=\Pi\big(s(g)s(h)\big)
\\
\!\!&=&\!\!
\Pi\big(z_{g,h}\,s(gh)\big)=\Pi\big(z_{g,h}\big)\Pi\big(s(gh)\big)=\Pi\big(z_{g,h}\big)\pi(gh), 
\end{eqnarray*}
and so, when $\Pi$ is irreducible, for $z\in Z$, we have $\Pi(z)=\chi_Z(z)I$, by {\it Schur's Lemma}, with $I$ the identity operator. We call the character $\chi_Z\in\widehat{Z}$ of $Z$ as {\it spin type}\, of $\Pi$ (with implicit reference to the base group $G\cong G'/Z$). Here 
$\Pi(z_{g,h})=r_{g,h}I$ with $r_{g,h}=\chi_Z(z_{g,h})\in{\boldsymbol C}^\times$, and accordingly 
\begin{eqnarray}
\label{2013-01-18-4}
\pi(g)\pi(h)=r_{g,h}\,\pi(gh)\quad(g,h\in G),
\end{eqnarray}
and so $\pi$ is a projective representation of $G\cong G'/Z$ with a factor set $r_{g,h}=\chi_Z(z_{g,h})$.

\vskip1em

{\bf 2.\; Represetation group.}\; 
A {\it representation group}\, $G'$ of a finite group $G$ is defined as a special central extension of $G$ such that any projective representation $\pi$ of $G$ can be obtained as above from a linear representation $\Pi$ of $G'$, and that the order $|G'|$ is minimum among central extensions having such property. Note that a projective representation of $G$ is nothing but a multi-valued representation. In the paper [Sch1], among other things, the following is shown: 
\\[1ex]
\indent
{\bf (1)}\; {\it  For an arbitrary finite group $G$, there exist finite number of representation groups (modulo isomorphisms),} 

{\bf (2)}\; {\it  For any representation group $G'$ of $G$, the commutative group $Z$ in the central extension (\ref{2013-01-19-2}) is isomorphic to Schur multiplier $H^2(G,{\boldsymbol C}^\times)$ of $G$. }
\bigskip

{\bf Note A.2.2.1.} In the case of a connected Lie group $G$, its universal covering groups correspond to the representation group in the case of finite groups. Taking a representation group of $G$, we denote it as $R(G)$. After the case of Lie groups, we call $R(G)$ as a {\it univrsal covering group}\, of $G$, even if it is not unique in general (Cf.\;{\bf A.2.3}, {\bf A.2.6}).

\vskip1em

{\bf 1907}\; {\bf （Construction of representation groups, their numbers, Schur multipliers and spin characters)}

\vskip.3em

{\small 
[Sch2]  
 J.\;Schur, {\it 
Untersuchungen \"uber die Darstellung der endlichen Gruppen  
 durch gebrochene lineare Substitutionen,}
J.\,f\"ur die reine und angew.\,Mathematik, {\bf 132}(1907), 85--137.
 }

 {\small \;\;[Title]\; Studies on the representations of finite groups by linear fractional transformations.}

 \vskip.6em
 
Here there were given methods of consturucting representation groups, evaluation of numbers of non-isomorphic representation groups, and a method of calulating Schur multipliers. Furthermore, irreducible spin characters (characters of spin irreducible representations) were explicitly given for groups 
 \\[1ex]
 \quad\hspace{5ex}
$SL(2,K),\,PSL(2,K),\,GL(2,K),\,PGL(2,K),$\, with $K=GF[p^n]$.   
\\[1ex]
\quad 
By this result, for these groups, the classification of irreducible spin representations were basically completed.

\vskip1.5em

{\bf 1911} {\bf （Construction of spin irreducible representations and calculations of spin irreducible characters for symmetric groups and alternating groups of $n$\,th order)} 

\vskip.2em
[Sch4]\;   J.\;Schur, {\it  
\"Uber die Darstellung der 
symmetrischen und der alternierenden Gruppen durch 
gebrochene lineare Substitutionen,} ibid., 
 {\bf 139}(1911), 155--255.

\vskip.3em
He constructed representation groups of the symmetric group 
$\mathfrak{S}_n$ of $n$-th order. Atually, for $\mathfrak{S}_n$ with $n=2,3$, its representation group is $\mathfrak{S}_n$ itself,\; and for $n\ge 4,\ne 6$, the number of non-isomorphic representation groups is 2，and for $n=6$ the number is 1. Also the type of representation groups of alternating group of $n$-th order is unique. 

Studying their linear representations, he gave spin representations of the base groups $\mathfrak{S}_n$ and $\mathfrak{A}_n$. Their construction and the calculation of their characters were explicitly carried out. 

For the construction of spin representations, he first gave a fundamental spin representation  $\Delta_n$ called {\it Hauptdarstellung}\, of $\mathfrak{S}_n$, utilizing almost the same $2\times 2$ matrices as the so-called Pauli matrix triplets (1927) in {\bf A.2.4} below. Well managing $\Delta_n$ as a seed, he succeeded to construct all the irreducible spin representations, and then calculated their characters.  

In the calculation of these spin characters, there appear the so-called $Q$-functions of Schur. In total, his results are inclusive and even today has important influence.

\vskip1.5em

\addcontentsline{toc}{subsection}{A.2.3\; \;Cases of Lie groups and Lie algegras} 

{\Large\bf \;A.2.3\; \;Cases of Lie groups and Lie algegras} 

\vskip.8em 

{\bf 1913} {\bf （Determination of irreducible linear representations of simple groups by means of the highest weights, and discovery of spin representations of rotation groups)}

\vskip.4em

{\small  
 [Car1]  \'E.\;Cartan, {\it 
Les groupes projectifs qui ne laissent invariante aucune multiplicit\'e plane,} 
  Bull. Soc. Math. France, {\bf 41}(1913), 53--96.
  
[Title] Groups of projective transformations which don't leave any linear subspaces invariant. 
}
 \vskip1em  
 
 \'{E}.\,Cartan\footnote{\'{E}lie Joseph Cartan, 9 April, 1869\;--\;6 May, 1951.} \!classifyed in this paper all the irreducible representation over ${\boldsymbol C}$ of any complex simple Lie algebras. Here he utilized his classification of simple Lie algebras $\mathfrak{g}$ over ${\boldsymbol C}$, and proved that an irreducible representation of $\mathfrak{g}$ is determined by its highest weight, and clarified how the highest weight is determined. 

At that time, since the results on the structues of Lie group $G$ corresponding to $\mathfrak{g}$ and its universal covering group are not sufficiently prepared, he did not so clearly recognize that, in the classified irreducible representations, there are contained spin representations. 
\vskip1em

For instance, the Lie algebra of 3-dimensional rotation group $SO(3)$ is 
 $\mathfrak{s}\mathfrak{o}(3)$, isomorphic to $\mathfrak{s}{\mathfrak u}(2)$ whose corresponding Lie group is $SU(2)$. The latter is a universal covering group of $SO(3)$. The 2-dimensional natural representation \,$SU(2)\ni u \mapsto u$\, is, if seen from $SO(3)$, a doubly-valued representation, and this is the origin of the name of {\it spin representation}. 

If we see this from the level of Lie algebras,  a complexification of $\mathfrak{s}\mathfrak{o}(3)$ is $\mathfrak{s}\mathfrak{o}(3,{\boldsymbol C})$, which is isomorphic to $\mathfrak{s}\mathfrak{l}(2,{\boldsymbol C})$.  A universal covering group correspondig to the latter is $SL(2,{\boldsymbol C})$, and its 
2-dimensional natural representation  \,$SL(2,{\boldsymbol C})\ni g\mapsto g$\, has the Lie algebra version  \,$\mathfrak{s}\mathfrak{l}(2,{\boldsymbol C})\ni X \mapsto X$\, and an appropriate basis of $\mathfrak{s}\mathfrak{l}(2,{\boldsymbol C})$ (triplet of matrices) is Pauli matirix triplet, explained below. 
 \vskip1em

 In the paper [Car1], regrettably enough, explicit matrix forms are not written, and so theoretical physicists (such as Pauli and Dirac) did not use the results in [Car1], and rediscover independently Pauli matrices again.

\vskip1em

{\bf Expository paper:\;} [CC]  
S.-S.\;Chern and C.\;Chevalley, 
{\it \'E.\,Cartan and his mathematical work,} 
  Bull. Amer.\;Math.\;Soc., {\bf 58}(1952), 217--250.

\vskip2em

\addcontentsline{toc}{subsection}{A.2.4\; \;Spin theory in quantum mechanics and mathematical foundation of quantum mechanics}

{\Large\bf \;A.2.4\; \;Spin theory in quantum mechanics and mathematical foundation of quantum mechanics}

\vskip.8em

{\bf 1927}\; {\bf （Discovery of Pauli matrix triplet and its applications)} 

\vskip.4em

[Pau] 
W.\;Pauli,\footnote{Wolfgang Ernst Pauli (25 April 1900\;--\;15 December 1958)} {\it 
Zur Quantenmechanik des magnetischen Elektrons,}  Zeitschrift f\"ur Physik,  {\bf 43}(1927), 601--623.

\vskip.5em

 Let $2\times 2$ Hermitian matrices be 
\begin{eqnarray}
\label{2011-01-17-1}
\qquad 
a=\sigma_1:= 
\begin{pmatrix} 
0 & 1 \\
1 & 0
\end{pmatrix}
, 
\quad 
b=\sigma_2:=
\begin{pmatrix} 
0 & -i \\
i & 0
\end{pmatrix}
, \quad
c=\sigma_3:= 
\begin{pmatrix} 
1 & 0 \\
0 & -1
\end{pmatrix}
, 
\end{eqnarray}
then these are {\it Pauli matrix triplet}\, (or {\it Pauli matrices}\, in short). Their commutation relations are 
\begin{eqnarray}
\label{2013-01-19-11}
\nonumber
\;
[a,b]=2ic,\quad [b,c]=2ia,\quad [c,a]=2ib\quad\;(i=\sqrt{-1}). 
\end{eqnarray}

\vskip.5em

{\bf A.2.4.1 \;Covering map $\Phi: SU(2) \to SO(3)$.}\; 

Triplet $\{a,b,c\}$ is a basis of $\mathfrak{s}\mathfrak{l}(2,{\boldsymbol C})$, isomorphic to the complexification of $\mathfrak{s}\mathfrak{o}(3)$. For more-mathematical expression, not containing $i=\sqrt{-1}$ in coefficients as above, take 
 \,$B_j:= i\sigma_j\;(j=1,2,3)$, then \,$\{B_1,B_2,B_3\}$\, is a basis over ${\boldsymbol R}$ of Lie algebra $\mathfrak{s}\mathfrak{u}(2)={\rm Lie}\big(SU(2)\big)$, and the commutation relations take the form 
\begin{eqnarray}
\label{2014-05-28-21}
 [B_j,B_k]=2B_l,
\end{eqnarray}
where $(j\;\,k\;\,l)$ is a cyclic permutation of $(1\;\,2\;\,3)$.
\vskip.5em 

Now, let us give some more mathematical explanation. As noted a little in Subsection {\bf A.2.3}, a covering map \,$\Phi:\,\widetilde{G}\ni u \mapsto g\in G$\, from the universal covering group $\widetilde{G}:=SU(2)$ to the rotation group $SO(3)$ is given as follows. For a column vector ${}^t(x_1,x_2,x_3)\in E^3$, we make correspond a $2\times 2$ Hermitian matrix as  
\begin{gather}
\label{2019-06-26-1}
X:=\sum_{1\leqslant j\leqslant 3}x_j\sigma_j= 
\begin{pmatrix}
x_3 & x_1-ix_2 
\\
x_1+ix_2& -x_3
\end{pmatrix}. 
\end{gather}

For $u\in SU(2)$, we give its action on $X$ as \,$X\mapsto X'=uXu^{-1}=uXu^*$ and express it as $X'=\sum_{1\leqslant j\leqslant 3}x'_j\sigma_j$. Then, we obtain a linear transformation as \,$(x_j)_{1\leqslant j\leqslant 3}\to (x'_j)_{1\leqslant j\leqslant 3}$. Its matrix expression gives the homomorphism  $g=\Phi(u)$. By simple calculations, we get the following formula (we leave the proof to readers):
$$
\Phi\big(\exp(\theta B_j)\big)=g_j(2\theta)\quad(1\le j\le 3),
$$  
where\; $\exp(\theta B_j),\,1\le j\le 3,$ and   
\;$g_j(\varphi),\,1\le j\le 3,$ are in this order 
\begin{gather*}
\label{2014-05-28-2}
\begin{pmatrix} 
\cos\theta & i\sin\theta
 \\
i\sin\theta & \cos\theta
\end{pmatrix}
, 
\quad\;\;
\begin{pmatrix} 
\cos\theta & \sin\theta
 \\
-\sin\theta & \cos\theta
\end{pmatrix}
, 
\quad\;\; 
\begin{pmatrix} 
e^{i\theta} & 0 
\\
0 & e^{-i\theta}
\end{pmatrix}
, 
\\[.5ex]     
\begin{pmatrix} 
1&0&0
\\
0&\cos\varphi & \sin\varphi 
\\
0&-\sin\varphi & \cos\varphi
\end{pmatrix}
, 
\;\;
\begin{pmatrix}
\cos\varphi&0&-\sin\varphi
\\ 
0&1&0
\\
\sin\varphi&0 & \cos\varphi
\end{pmatrix}
, 
\;\; 
\begin{pmatrix} 
\cos\varphi&-\sin\varphi&0
\\
\sin\varphi&\cos\varphi & 0 \\
0&0 & 1
\end{pmatrix}
. 
\end{gather*}

In this case, ${\rm Ker}(\Phi)=\{\pm E_2\}$ and so $\Phi$ gives a double covering. Moreover, $g_j(\varphi)$ is the matrix of a rotation arround the central axis $x_j$ of angle $\varphi$, for $j=1,2$, that advances the right-handed screw in the positive direction of the central axis, and for $j=3$, that advances the left-handed screw. Concerning one-parameter subgroup $\exp(\theta B_j)$, if we take one cycle $0\le \theta \le 2\pi$ of angle $\theta$, then its image $g_j(\varphi),\, \varphi=2\theta,$ rotate two cycles (arround the axis $x_j$). The generators of one-parameter subgroup $g_j(\varphi)$ of $G$ are  
\begin{gather}
\label{2019-06-30-1}
A_j:=\frac{d}{d\varphi}g_j(\varphi)\big|_{\varphi=0}\;\in\;\mathfrak{so}(3)=\mathfrak {so}(3,{\boldsymbol R})\quad(j=1,2,3).
\end{gather} 

The set $\{A_1,A_2,A_3\}$ gives a basis of $\mathfrak{so}(3)$, and has a system of commutation relations as  \,$[A_j,A_k]=A_l$, where $(j,k,l)$ is a cyclic permutation of $(1,2,3)$. The differential $d\Phi$ gives an isomorphism  \,$\tfrac{1}{2}B_j\to A_j\;(1\le j\le 3)$\, from $\mathfrak{m}\mathfrak{u}(2)\to \mathfrak{s}\mathfrak{o}(3)$.   
 
\vskip1.2em

{\bf A.2.4.2 \;Two-dimensional spin representation of $SO(3)$.}\; 

For $g\in G=SO(3)$, take a preimage $u\in \widetilde{G}=SU(2),\,\Phi(u)=g,$ of $\Phi$ and put 
\begin{gather}
\label{2019-06-30-2}
\pi_2(g):=u,
\end{gather}
then this is the 2-dimensional spin representation of the rotation group $G$. 
In reality, if $g$ moves in a small nieghbourhood of the identity $E_3$, then $g\mapsto \pi_2(g)$ gives a  1\,:\,1 map (local isomorphism from $G$ to $\widetilde{G}$), but if $\varphi$ moves continuously from $\theta_0$ to $\varphi=\theta_0+2\pi$, then it arives to \,$
\pi_2\big(g_j(\theta_0+2\pi)\big)=-\pi_2\big(g_j(\theta_0)\big)
$, 
that is, $-1$ appears as a multiplication factor. Thus it is natural to accept this as a doubly-valued representation ({\it spin}\, representation), or it can only be regarded as such. 

However when we step up from $G$ to its universal covering group $\widetilde{G}$, it becomes univalent linear representation (the identity representation $u\mapsto u$). 

A maximal commutative subgroup (one of Cartan subgroups) of $G=SO(3)$  is given by 
\begin{gather}
\label{2019-06/30-3}
H=\{g_3(\varphi)\,;\,0\le \varphi<2\pi\}.
\end{gather}

For a representation $\pi$ on a representation space $V(\pi)$, take a common eigenvector ${\boldsymbol v}\ne 0$ for the commutative family $\pi(h)\;(h\in H)$, and call it as a {\it weight vector}\, of $\pi$ with {\it weight}\, $\chi\in\widehat{H}$, where  $\pi(h){\boldsymbol v}=\chi(h){\boldsymbol v}\;(h\in H)$.　

For 
$h=g_3(\varphi)$, there exists a $k\in \tfrac{1}{2}{\boldsymbol Z}$ such that $\chi(h)=e^{ik\varphi}$. Suppose $\pi$ is irreducible. For $\pi$ (one-valued) linear representation, $k$ is integral $(k \in{\boldsymbol Z})$, and for (2-valued) spin representation, $k$ is half an odd integer. The maximal one among these weights $k$ is called the {\it highest weight}\, of $\pi$. We can prove by simple calculations that the highest weight determines the equivalence class of $\pi$ completely. 

Note that the highest weight of the above two-dimensional $\pi_2$ is equal to $1/2$.

\vskip1.5em

{\bf A.2.4.3 \;Original reason for the discovery of Pauli matrices.}\; 
\vskip.3em

The reason why Pauli arrived to the results in the above paper is the following:\begin{quotation}
\noindent 
In the principle that 
\lq\lq\,Catch an electron as a mass point\,'', the theory of applying wave equation to a wave function gives double number as the number of stable orbits of an electron arround an atomic nucleus, and this phenomena is called {\it duplexity phenomena}.  To avoid this contradiction, Pauli introduced the fourth quantum number, saying that electron has spin angular momentum of $\pm\tfrac{1}{2}\hbar$. 
\end{quotation}

To say this mathmatically, 
\lq\lq\,Wave equation and wave function of electron receive transformations not only by rotations in the space (as coordinate transformations) under the group of rotations $G=SO(3)$ but also transformations (through spin representation $\pi_2$) under its universal covering group $\widetilde{G}=SU(2)$. So it means that 
\begin{quotation}
{\it 
\noindent
electron does not live in the 3-dimensional Eucledean space but in the 
 2-dimensional space ${\boldsymbol C}^2$ (or 4-dimensional space over ${\boldsymbol R}$) on which $\widetilde{G}$ acts.
}
\end{quotation}

This too much purely-mathmatical explanation could not be easily accepted by physisists, and for instance, one of my university class mates who studied physics said that, in the lecture course of physics in Departement of Physics, Kyoto University, Professor expalined in such a way that, {\it arround some axis (of which direction we don't know), electron rotates}.  

Detailed description and explanations about Pauli's wave equation and wave functions and how is the transformation of the group $\widetilde{G}=SU(2)$, are given in Chap.\;4 of our book [HiYa]. 

\vskip1.5em

{\bf 1928}\; {\bf （Quantum Theory of Electron）} 
\vskip.2em
{\small 
[Dir] 
P.M.A.\;Dirac, {\it 
The quantum theory of electron,}  
 Proceedings of the Royal Society London A, {\bf 117}(1928), 
 610-624; and\; Part II,   ibid., {\bf 118}(1928),   351--361. 
 } 
\vskip.3em

In this paper, Dirac enlarged a system of three matrices of Pauli to a system of four matrices and obtained spin representation of Lorentz group, and discovered the so-called Dirac equation which is a relativistic (that is, Lorentz-invariant) wave equation. When electron is expressed with this equation the above mentioned duplexy phenomena disappear automatically without any assumption. 

Moreover Klein-Gordon equation is factorized into two factors of \lq Dirac equations', and the latter is known to be much fundamental. Dirac equation gives a foundation of later big development of quantum mechanics. To explain this in more detail, we first prepair some fundamentals about {\it Minkowski space}\, $M^4$ and {\it Lorentz group}\, ${\cal L}_4=SO_0(3,1)$.

 \vskip1.2em
 
 {\bf A.2.4.4 \;Minkowski space $M^4$.}\; 

The Minkowski space $M^4$ is a space to describe (space, time) consisting of (column) vectors \;${\boldsymbol x}={}^t(x_1,x_2,x_3,x_4)$ with 
 \begin{gather*}
 {}^t(x_1,x_2,x_3)\in E^3,\;\;x_4=ct\in {\boldsymbol R},\;c=\mbox{\rm light velocity},\;\;t=\mbox{\rm time},  
\end{gather*}
and equipped with inner product given as 
 \begin{gather*}
\label{2011-02-10-11}
\langle \mbox{\boldmath$x$},\mbox{\boldmath$x$}\rangle_{3,1}:=x_1^{\;2}+x_2^{\;2}+x_3^{\;2}-x_4^{\;2}={}^t{\boldsymbol x}\,J_{3,1}{\boldsymbol x}\,, 
\quad
 J_{3,1}:=\mathrm{diag}(1,1,1,-1). 
\end{gather*}
The group $O(3,1)$ consists of all $g \in GL(4,{\boldsymbol R})$ leaving the inner product invariant:  
$$
\langle g\mbox{\boldmath$x$},g\mbox{\boldmath$x$}\rangle_{3,1} = \langle \mbox{\boldmath$x$},\mbox{\boldmath$x$}\rangle_{3,1}\quad(\mbox{\boldmath$x$}\in M^4),$$
or\; ${}^tg J_{3,1}g=J_{3,1}$. It has 4 connected components and the one containing the identity element $e=E_4$ is called as {\it proper Lorentz group}\, and given as  
 \begin{eqnarray}
\label{2011-02-10-12}
\nonumber
&& 
{\cal L}_4:=SO_0(3,1):=\big\{g=(g_{ij})_{1\leqslant i,j\leqslant 4} \in O(3,1)\;;\; \det(g)=1,\;g_{44}\ge 1\big\}.
\end{eqnarray}

\vskip.3em

{\bf A.2.4.5 \;Dirac equation.}\; 

Dirac's wave function $\psi(\mbox{\boldmath$x$})$ is a function on $M^4$ with values in $V={\boldsymbol C}^4$, written as column vector as \;$\psi(\mbox{\boldmath$x$})=\big(\psi_j(\mbox{\boldmath$x$})\big)_{1\leqslant j\leqslant 4}$. Dirac gave the so-called {\it Dirac equation}\, for elelctron, {\it relativistic}, as follows: put 
  \begin{eqnarray}
\label{spin30}
&&
\sigma_1 = 
\left(\!
\begin{array}{rr}
0 & 1 \\
1 & 0
\end{array}
\!
\right), \,
\sigma_2 = 
\begin{pmatrix}
0 & -i \\
i & 0
\end{pmatrix}
, \,
\sigma_3 = 
\begin{pmatrix}
1 & 0 \\
0 & -1
\end{pmatrix}
, \,
\sigma_4 = 
\left(\!
\begin{array}{rr}
1 & 0 \\
0 & 1
\end{array}
\!
\right), 
\end{eqnarray}
and further put $4\times 4$ matrices\; 
  $\gamma_1,\gamma_2,\gamma_3,\gamma_4$ as 
\begin{eqnarray}
\label{dirac3}
&&
\gamma_j = 
\begin{pmatrix}
0_2 & -i\sigma_j \\
i\sigma_j & 0_2
\end{pmatrix}
\;\;\;(1 \le j \le 3),\quad 
\gamma_4 = 
\begin{pmatrix}
-i\sigma_4 & 0_2 \\
 0_2 & i\sigma_4
\end{pmatrix}
,
\end{eqnarray}
where $0_2$ denotes zero matrix of order 2. Then, in the case where the electromagnetic field is zero, Dirac equation is given, with the mass $m$ of electron, as 
\begin{eqnarray}
\label{2011-02-11-32}
\qquad
\big(D+\kappa\big)\psi=0,\quad 
D:=\!
\sum_{1\leqslant j\leqslant 4}\!\gamma_j\partial_{x_j},\;\;\partial_{x_j}:=\frac{\partial\;}{\partial x_j},\;\;\kappa:=\frac{mc}{\hbar}\,.
\end{eqnarray}

The action of $g\in{\cal L}_4$ on the wave function is  \;$T_g(\psi)({\boldsymbol x}):=g\big(\psi(g^{-1}{\boldsymbol x})\big)$. Dirac equation is {\it Lorentz-invariant}\, or ${\cal L}_4$-invariant, which is expressed as 
$$
T_g^{\;-1}\big(D+\kappa\big)\big(T_g\psi\big)=\big(D+\kappa\big)\psi\quad (g\in{\cal L}_4).
$$

Moreover, Klein-Gordon equation is factorized as  
\begin{eqnarray*}
\quad
(\square -\kappa^2)\psi=(D-\kappa)(D+\kappa)\psi=0, \qquad 
\square :=\! \sum_{1\leqslant j \leqslant 3}\!\partial_{x_j}^{\;\,2}-
\partial_{x_4}^{\;\,2}\,.
\end{eqnarray*}

For more detailed explanation in this direction, please cf.\;[Ch.\,5, HiYa] and [HHoH].

\vskip1.5em

\noindent

{\bf 1928} \;{\bf (Mathematical foundation of quantum mechanichs (Neumann),\\ 
\hspace{10ex} Group actions are projective representations (Weyl)\,)}

\vskip.8em

{\small

[Wey2]\; H.\;Weyl,
{\it Gruppentheorie und Quantenmechanik}, 1st edition, 1928\,; 
 (2nd edition, 1931), Hirzel, Leipzig.\quad English transl. of 2nd ed. by H.\,P. Robertson, Dutton, N.Y., 1932. 
\vskip.4em

 [Neu] J.\,von Neumann, {\it Mathematische Grundlagen der Quantenmechanik}, J.\,Springer, 1932. 

English transl. by R.\,T. Beyer: {\it Mathematical Foundations of Quantum Mechanics}, Princeton University Press, 1955.

}

 \vskip.8em
Following von Neumann's idea (see [Neu], for instance) to describe the state of an elementary particle, one utilizes a unit vector in a Hilbert space, or more exactly a complex line (Strahlenk\"orper) determined by the vector. Accordingly if a group acts on the space of states, then the action of this transformtion group is described by its {\bf Strahldarstellung\;=\;ray representation}, which is essentially equal to a {\it projective} or {\it spin representation.}

\vskip2em

\addcontentsline{toc}{subsection}{A.2.5\; \;Rediscovery of a work of A.H.\;Clifford}

{\Large\bf \;A.2.5\; \;Rediscovery of a work of A.H.\;Clifford}
\vskip.8em

\thispagestyle{plain}

{\bf 1937} {\bf （Situation when an irreducible representation of a finite group is\\ 
\hspace*{6.7ex} restricted to its normal subgroup）} 
\vskip.4em

{\small 
[Clif] A.H.\;Clifford, {\it 
Representations induced in an invariant subgroup,} 
Ann. Math., {\bf 38}(1937), 533--550. 
 }

\vskip.4em

A.H.\;Clifford is an american mathematician, and his life time don't overlap with that of well-known W.K.\;Clifford of Clifford algebra. He was an assistant of H.\;Weyl, at Princeton Institute for Advanced Study, when he wrote this paper. In it he studied in detail the restriction onto a normal subgroup $N$, of an irreducible representation of a finite group $H$. He found, in particular, in some cases, there appear inevitably {\it projective} (or {\it spin}) representations as tensor product factors.  

The following theorem is a result of summarizing some essential parts in \S\S 3--4 of this paper in my own way.  

\vskip1.2em 

{\bf Theorem A.2.5.1.} (cf. [Clif, \S\S 3--4])\;\; 
{\it Let  
$N$ be a normal subgroup of a finite group $H$. For an irreducible linear representation  $\tau$ of $H$, suppose that all the irreducible components of the restriction $\tau|_N$ are mutually equivalent, that is, \;$\tau|_N\cong [\ell]\cdot \rho$,where $\rho$ is an irreducible representation of $N$, and  
\;$\ell=\dim \tau/\dim \rho$\, is its multiplicity. Moreover assume that the base field is algebraically closed.  

{\rm (i)}\; Irreducible representation 
 $\tau$ of $H$ is equivalent to a tensor product of two irreducible projective matrix representations\, $C$ and\, $\Gamma$\,:
\begin{equation}
\label{2010-03-14-8} 
\tau(h) \cong \Gamma(h)\otimes C(h)\quad(h\in H),   
\quad 
\end{equation}
where the factor sets of\, $C$ and $\Gamma$ are mutually inverse of the other, and\\
\hspace{17ex}
$\dim \Gamma=\ell,\quad \dim C=\dim \rho,$  
\vskip.3em
{\qquad\qquad} 
$\rho(h^{-1}uh)=C(h)^{-1}\rho(u)C(h)\quad(h\in H,\,u\in N)$. 

\vskip.8em

{\rm (ii)}\; It can be normalized as follows: for\, 
$h\in H,\;u\in N$, 
\begin{eqnarray*}
C(hu) =C(h)\rho(u),\;C(u)=\rho(u)\,;\quad \Gamma(hu)=\Gamma(h),\;\;\Gamma(u)=E_\ell.
\end{eqnarray*}
In this situation, $h\mapsto \Gamma(h)$ is substantially a projective representation of the quotient group $H/N$, and the factor sets of\, $C$ and\, $\Gamma$ are mutually inverse of the other, and they are substantially factor sets of $H/N$.

}
\vskip.3em 

\bigskip

The result of Clifford has an important meaning for the stand point of spin representations in the theory of group representations. It says that evenwhen we are working only on {\it linear} representations of groups, there may appear inevitably {\it spin} representations of some subgroups in the process of constructing all irreducible linear representations.

One simple such example is the case of semidirect product groups $H=U\rtimes S$ in [Hir1], where $U$ is compact, normal in $H$, and $S$ is finite. In [Hir1] we gave a method of constructing a set of complete representatives of the dual $\widehat{H}$ by using {\it spin} representations of certan subgroups of $S$. Our proof there is independent of the result of Clifford but we can give another proof of it by utilizing Theorem A.2.5.1 above. 

Note that, in the work of G.W.\;Mackey\,[Mac1]--[Mac2], $S$ is assumed to be abelian and so there do not appear any {\it spin} representations. The above classical results of Clifford is also linked to our recent works [HHH4], [HHo] and [Hir2].

\vskip1.2em

{\bf Note A.2.5.1.}\; Now let me give a short break, intermediate. Then I would like to state some personal thought about {\it spin representations}. By an invitation of late Prof.\;Mitsuo Sugiura, I joined to series of anual meetings held at Tsuda University as \lq\lq\,Symposium on History of Mathematics\,''. In the first year, I reported on my work on representations and characters of semisimple Lie groups as a part of mathematics at present. But starting from the next year, I reported on a series of meomorial works of Frobenius which originated the theory of linear representations and chracters of (finite) groups, according to the years of publication successively as in [H1], consecutively in four years. 

After that, I begun to follow works of Schur, a student of Frobenius. Papers of his first stage are some kind of extention or simplifiction of the results of his teacher. Then, at the second stage, I met the difficult work, a triplet of papers ({\small [Sch1, 1904], [Sch2, 1907], [Sch4, 19011])}, on {\it spin} (or {\it projective} or {\it multi-valued}) representations of groups, in the terminology of Schur, {\it Darstellung durch gebrochene lineare Substitutionen}, and I made 2 reports on these papers in [H2]. 

Then my interest spreaded naturally over the surrounding area containing the case of Lie groups, and I felt curious about something. In the 3rd paper [Sch4] of Schur,  there appeared fully, many interesting results on spin (or projective) representations of symmetric groups  ${\mathfrak S}_n$ and alternating groups ${\mathfrak A}_n$. However, after that, there is a significant long, blank period before some succeeding results were published on {\it finite groups}, whereas, for {\it semisimple Lie groups}\, including motion groups and Lorentz groups, the theory of spin representations was studied continuously and steadily. 

For instance, after almost half a century of blank, there finally appear Morris' paper [Mor1, 1962] on spin representations of symmetric groups, and Ihara-Yokonuma's paper [IhYo, 1965] on Schur multipliers ${\mathfrak M}(G)$ for finite (and also infinite) reflexion groups $G$. For these historical phenomena, there might be some reason, but I felt that the theory of spin representations in this area are treated in some sense as a stepchild of the main theory of representations. 

However, in the above work of A.H.\,Clifford, it is proved that, in the theory of linear representations of $G$, there appears naturally and inevitably {\it spin} (or {\it multi-valued} or {\it projective}) representations of some subgroups of $G$. This shows that {\it spin} representations are nothing but a legitimate child of the theory of group representations.

\vskip2em

\addcontentsline{toc}{subsection}{A.2.6\; \;Development in mathematics, and quantum mechanics}

{\Large\bf \;A.2.6\; \;Developments in Mathematics, and Quantum\\ 
 \hspace*{9.7ex} Mechanics}

\vskip.8em

{\bf 1939}\; {\bf （Construction of all irreducible representations and calculation of their characters)}

\vskip.3em
{\small 
[Wey1] H.\;Weyl,
{\it Classical Groups, Their Invariants and 
Representations}, Princeton University Press, 1939. 
}

\vskip.4em

Classical groups (or Lie groups of classical type) over ${\boldsymbol C}$ are given, up to isomorphisms, as follows: 
 \\
${\qquad\qquad} 
{\rm type}\;\;A_n\qquad\;SL(n\!+\!1,{\boldsymbol C})\;\;\quad\;(n\ge 1), 
\\
{\qquad\qquad} 
{\rm type}\;\;B_n\qquad\;SO(2n\!+\!1,{\boldsymbol C})\;\;\hspace{1.4ex}(n\ge 2),
\\
{\qquad\qquad} 
{\rm type}\;\;C_n\qquad\;Sp(2n,{\boldsymbol C})\;\;\hspace{5.5ex}(n\ge 3),
\\
{\qquad\qquad} 
{\rm type}\;\;D_n\qquad\;SO(2n,{\boldsymbol C})\;\;\hspace{4.5ex}(n\ge 4),
$ 
\\[.3ex]
and their compact real forms are respectively 
\begin{eqnarray*}
\;
SU(n\!+\!1),\;\;SO(2n\!+\!1),\;\;Sp(2n)=SU(2n)\cap Sp(2n,{\boldsymbol C}),\;\;SO(2n).
\end{eqnarray*}

By means of the so-called {\it Weyl's unitarian trick}, irreducible unitary representations of the latter corresponds 1-1 way to irreducible (finite-dimensional) holomorphic representations of the former. These representations are all constructed at least in principle and their characters are explicitly calculated.

The compact groups $SU(n\!+\!1)$ and $Sp(2n)$ are simply connected together with their mother groups $SL(n\!+\!1,{\boldsymbol C})$ and $Sp(2n,{\boldsymbol C})$. Whereas $SO(N),\,N\ge 3,$ has (double covering) universal covering group ${\rm Spin}(N)$ and so they have irreducible linear representations and also irreducible {\it spin} (or {\it projective}, or {\it doubly-valued}) representations. But they can be treated in a same way, and have a unified common formula for characters and dimensions for each type. 
\vskip1.2em

{\bf 1947}\; 
{\bf (Irreducible representations and charactrs of $SL(2,{\boldsymbol R})$, double covering of 3-dimensional Lorentz group)}

\vskip.3em
{\small 
[Bar]  
V.\;Bargmann, {\it  
Irreducible unitary representations of the Lorentz group,}
 Ann. Math., {\bf 48}(1947), 568--640. 
 }
 \vskip.4em 

A double covering group of 3-dimensional Lorentz group ${\cal L}_3=SO_0(2,1)$, with the 2-dimensional space-part, is realized by $SL(2,{\boldsymbol R})$. Bargmann constructed all irreducible unitary (and also non-unitary) representations $SL(2,{\boldsymbol R})$ and gave their characters. Seeing from the base group ${\cal L}_3$, a half of these representations is linear and another half is {\it of spin} (doubly-valued). A universal covering group of ${\cal L}_3$ is $\infty$-times covering and the situation differs very much from the case of $SL(2,{\boldsymbol R})$.

 \vskip1.2em

{\bf 1947}\; 
{\bf (Irreducible representations of $SL(2,{\boldsymbol C})$, a universal covering group of 4-dimensional Lorentz group)}

\vskip.3em
{\small [GeNa] I.M.\;Gelfand and M.I.\;Naimark, {\it 
Unitary representations of Lorentz
 group}\, (in Russian), Izvestia Akad. Nauk SSSR, {\bf 11}(1947), 411--504 [\,English Translation in Collected Works of Gelfand, Vol.\;2, pp.41--123\,]. 
 } 

\vskip.4em
The universal covering group $SL(2,{\boldsymbol C})$ of 4-dimensional Lorentz group ${\cal L}_4$ is doubly covered. In this paper, all irreducible representations are constructed and their characters are calculated explicitly. Seeing from the base group ${\cal L}_4=SO_0(3,1)$, there are linear representations and {\it spin}\, representations, and they can be treated in a similar way.

\vskip1.2em

{\bf 1947}\; 
{\bf (Infinitesimal construction of irreducible representations of 4-dimensional Lorentz group)} 

\vskip.3em
{\small 
[HC]  
Harish-Chandra (Doctor's Thesis), {\it 
Infinite irreducible representations of the Lorentz 
 group,}  Proceedings of the Royal Society London A, {\bf 189}(1947), 327--401.
  }

 \vskip.4em 

In colonial ages of India, an excellent local academic person could be sent to colonial mother country, Great Britain, for much advanced studies. Thus Harish-Chandra studied in London for his Doctrate under the supervision of Dirac. 

Anyhow, the subject of this work was suggested to him by Dirac as is noted in its Introduction. 
He classified and constructed irreducible representations of 4-dimensional Lorentz group ${\cal L}_4=SO_0(3,1)$ on the level of Lie algebras. This is called as {\it infinitesimal}\, method. Since this is the story on the level of a Lie algebra, there didn't appear the theory of characters, which are special functions on the corresponding Lie group.
 
 \vskip2em

\addcontentsline{toc}{subsection}{A.2.7\; \;Renaissance of the theory of spin representations}

{\Large\bf \;A.2.7\; \;Renaissance of the theory of spin representations of finite groups}

\vskip.8em 

 For a connected semisimple Lie group, when we treat finite dimensional irreducible representations, their highest weights have some differences between the cases of linear representations and spin (multi-valued) representations, but in principle they can be treated equally in a unified way. 

For example, for the rotation group $SO(3)$, the highest weights of the former are non-negative integers, whereas those for doubly valued spin representations are positive half integers. However the character formula and the dimension formula take the same form and have no difference. And passing through the ages of \'{E}.\;Cartan, J.\;Schur (= I.\;Schur), and H.\;Weyl, of finite-dimensional representations, 
we came to such ages as treating infinite-dimensional representations of Lorentz groups and their covering groups  $SL(2,{\boldsymbol R})$ and $SL(2,{\boldsymbol C})$ in a natural unified way, not so much depend on {\it linear}\, or {\it spin}, and so on. 

On the contrary, for finite groups, linear representations and {\it spin}\, (or {\it multi-valued}) ones are, in many cases, have quite different features. Maybe depending on this fact partly, and also depending partly on that Schur worked out in his trilogy  [Sch1], [Sch2], and [Sch4], too completely, comprehensively and thoroughly, so that his successors could not appear for long time.

\vskip1.2em

{\bf 1962}\; 
{\bf (Beginning of recapturing Schur's theory)}

\vskip.3em
{\small 
[Mor1] A.O.\;Morris, {\it 
The spin representation of the symmetric group,} Proc.   
 London Math. Soc., {\bf (3)\;12}(1962), 55--76. 
 }
\vskip.4em

After half a century of Schur's trilogy, there opened the study of spin representations. The above is the first of such papers. Its content begun with recapturing Schur's theory. Later, students of Morris and other peoples gradually follow spin theory of finite groups. 

\vskip1.2em

{\bf 1965}\; 
{\bf (Calculation of Schur multipliers $H^2(G,{\boldsymbol C}^\times)$)} 

\vskip.3em
{\small 
[IhYo] 
  S.\;Ihara and T.\;Yokonuma, {\it 
 On the second cohomology groups (Schur
 multipliers) of finite reflexion groups,} 
  J.\,Fac.\,Sci.\,Univ.\,Tokyo, Ser.1,   {\bf IX}(1965), 155--171. 
}
\vskip.4em

From a more algebraic point of view, the determination of Schur multiplier $H^2(G,{\boldsymbol C}^\times)$ is the theme studied hardly and intensively in some period of time. According to Gorenstein's book [Chap.2,\, Gor], Schur multiplier is treated as an important datum for known finite groups to prepare a complete classification of finite simple groups, and this gives partly a motivation of calculating $H^2(G,{\boldsymbol C}^\times)$ explicitly. 

However, from the explicit determination of $H^2(G,{\boldsymbol C}^\times)$ to the explicit determination of all irreducible {\it spin} (multi-valued) representations, they didn't proceed immediately to this second step.  There was certain time gap. Why\;?\,  I don't know, but is there some psychological barrier ?\!
\footnote{
Once I wrote a letter to Prof.\,Yokonuma thanking his two papers on Schur multipliers of finite and infinite reflexion groups. At that occasion I put a similar question like \lq\lq\,Why don't you proceed, after determination of Schur multipliers, to construction process of multi-valued irreducible representations, taking the datum $H^2(G,{\boldsymbol C}^\times)$ for help ?\,'' His answer is somewhat vague and means something like {\it \lq\lq\,If I have such an idea, then it would be much better for me.''}}  

Here I quote early three papers. 

\medskip

{\bf 1973}\; 
{\small 
[Mor2]  A.O.\;Morris, {\it 
Projective representations of abelian groups,}  
 London Math. Soc., {\bf (2) 7}(1973), 235--238.
}
\medskip

{\bf 1974}\; 
{\small [DaMo] 
J.W.\;Davies and A.O.\;Morris, {\it 
The Schur multiplier of the generalized symmetric group,}  J. London Math. Soc., {\bf (2) 8}(1974), 615--620. 
}
\medskip

{\bf 1976}\; 
{\small 
[Rea]  
E.W.\;Read, {\it 
On the Schur multipliers of the finite imprimitive 
 unitary reflexion groups $G(m,p,n)$,} 
  J. London Math. Soc., {\bf (2), 13}(1976), 
 150--154. 
 }

 \vskip2em

\addcontentsline{toc}{subsection}{A.2.8\; \;Weil representations of symplectic groups}

{\Large\bf A.2.8\; \;Weil representations of symplectic groups} 
\vskip.8em 

{\bf 1964}\; 
{\bf (Weil representations of symplectic groups over locally compact abeian groups)}

\vskip.3em
{\small 
[Wei2] 
A.\;Weil, {\it 
Sur certains groupes d'op\'erateurs unitaires,} 
 Acta Math., {\bf 111}(1964), 143--211. 
 }

\medskip

As for the motivation of this work, maybe Weil had as his purpose cetain images of working on algebraic groups over adele groups. Here, more generally, he uses irreducible representations of 
Heisenberg groups over locally compact abelian groups and groups consisting of their intertwining operators. These groups give spin representations of symplectic groups, which is called {\it Weil representations}.  

In the special case where the locally compact abelian group is simply real number field ${\boldsymbol R}$, the symplectic group appeared is just $Sp(2n,{\boldsymbol R})$ and the group of intertwining operators is its double covering group, called {\it Metaplectic group}\, and denoted as $M\!p(2n,{\boldsymbol R})$. 
Actually, Weil representation in this case is a doubly-valued representation of the base group $Sp(2n,{\boldsymbol R})$, which turns out to be a linear representation of $M\!p(2n,{\boldsymbol R})$ if going up to the upper level.

\vskip1.2em

{\bf 1972}\; 
{\small [Sait]  
M.\;Saito, {\it 
Repr\'esentations unitaires des groupes symplectiques,} 
J. Math. Soc. Japan, {\bf 24}(1972), 232--251. 
 } 
\vskip.3em

A maximal compact subgroup of $Sp(2n,{\boldsymbol R})$ is isomorphic to $U(n)$. Their universal covering groups $\widetilde{S\!p}(n,{\boldsymbol R})$ and $\widetilde{U}(n)$ are both infinite times covering. For any positive integer $k$, there exists (exactly) $k$-times covering group. 
Weil representation is doubly-valued for $Sp(2n,{\boldsymbol R})$, and linear for $M\!p(2n,{\boldsymbol R})$. M.\;Saito constructed several series of irreducible unitary representations, using Weil representation\,:

\vskip1.8em

$
\qquad
\hspace{17ex} 
\begin{array}{l}
\;\widetilde{Sp}(2n,{\boldsymbol R})\quad \mbox{\rm \;the universal covering group}
\\
\;\;\;\quad\, 
\downarrow
\\
\,M\!p(2n,{\boldsymbol R})\quad \mbox{\rm \;double covering group}
\\
\qquad
\downarrow
\\
\;Sp(2n,{\boldsymbol R})\quad \mbox{\rm \;base group}
\end{array}
$
\\[2ex]

{\bf 1979}\; 
{\small 
[Yos1] H.\;Yoshida, {\it 
Weil's representations of the symplectic groups over 
 finite fields,} J. Math. Soc. Japan, {\bf 31}(1979), 399--426.
 }
 \vskip.4em 
 
 For the symplectic group $Sp(2n,K)$ over a finite field $K$, its Weil representation is an irreducible linear representations of $Mp(2n,K)$.

\vskip1.2em

{\bf 1992}\; 
{\small 
[Yos2]  H.\;Yoshida, {\it 
Remarks on metaplectic representations of $SL(2)$,} 
J. Math. Soc. Japan, {\bf 44}(1992), 351--373.
}
\vskip.4em

For a local field with characteristic $\ne 2$, consider $n$-times covering group $\widetilde{G}$ of $G=SL(2,K)$. For $n=2$ and $K\ne {\boldsymbol C}$, Weil representation of $G$ can be lifted up to an irreducible linear representation of $\widetilde{G}$. In this paper, for any $n\ge 3$, Yoshida constructed spin irreducible representations of $G$ which becomes linear if going up to exactly $n$-times covering $\widetilde{G}$.

{\small

\end{document}